\documentclass[12pt]{amsart}
\usepackage{amsfonts,amssymb,graphicx}

\newtheorem{theorem}{Theorem}

\newtheorem{corollary}[theorem]{Corollary}

\theoremstyle{definition}
\newtheorem{definition}[theorem]{Definition}
\newtheorem{example}[theorem]{Example}
\newtheorem{notrems}[theorem]{Notation and Remarks}

\newtheorem{remark}[theorem]{Remark}

\newcommand{\Section}[1]{\section{#1}\setcounter{theorem}{0}}

\newcommand{\<}{\langle}
\renewcommand{\>}{\rangle}
\newcommand{\SO}{{\mathit{SO}}}
\renewcommand{\O}{{\mathcal O}}
\renewcommand{\S}{{\mathcal S}}
\newcommand{\V}{{\mathcal V}}
\renewcommand{\Gamma}{{\varGamma}}
\newcommand{\R}{{\mathbb R}}
\newcommand{\Z}{{\mathbb Z}}
\newcommand{\N}{{\mathbb N}}
\newcommand{\Id}{{\operatorname{Id}}}
\newcommand{\Iso}{{\operatorname{Iso}}}
\newcommand{\tr}{{\operatorname{tr}}}
\newcommand{\dimm}{{\operatorname{dim}}}
\newcommand{\spann}{{\operatorname{span}}}
\newcommand{\trans}[1]{{\,}^t\hskip-1.3pt#1}
\newcommand{\inv}{^{-1}}

\begin{document}

\title[Isospectral orbifolds with different maximal isotropy orders]
{Isospectral orbifolds with different maximal isotropy orders}


\author{Juan Pablo Rossetti}
\address{Famaf-CIEM, Universidad Nacional de C\'ordoba, 5000 C\'ordoba, Argentina}
\email{rossetti@mate.uncor.edu}
\author{Dorothee Schueth}
\address{Institut f\"ur Mathematik, Humboldt-Universit\"at zu
Berlin, D-10099 Berlin, Germany}
\email{schueth@math.hu-berlin.de}
\author{Martin Weilandt}
\address{Institut f\"ur Mathematik, Humboldt-Universit\"at zu
Berlin, D-10099 Berlin, Germany}
\email{weilandt@math.hu-berlin.de}


\thanks{All three authors were partially supported by DFG Sonderforschungsbereich~647.}

\begin{abstract}
We construct pairs of compact Riemannian orbifolds which are isospectral for the
Laplace operator on functions such that
the maximal isotropy order of singular points in one of the orbifolds
is higher than in the other. In one type of examples, isospectrality
arises from a version of the famous Sunada theorem which also implies
isospectrality on $p$-forms; here the orbifolds
are quotients of certain compact normal homogeneous spaces.
In another type of examples, the orbifolds are quotients of Euclidean $\R^3$
and are shown to be isospectral on functions using dimension formulas
for the eigenspaces developed in \cite{MR:2001}.
In the latter type of examples the orbifolds are not isospectral on $1$-forms.
Along the way we also give several additional examples of isospectral
orbifolds which do not have maximal isotropy groups of different size
but other interesting properties.
\end{abstract}

\maketitle

\Section{Introduction}
\label{sec:intro}

\noindent
This paper is concerned with the spectral geometry of compact Riemannian orbifolds.
The notion of Riemannian orbifolds is a generalization of the notion of Riemannian
manifolds. In a Riemannian orbifold each point has a neighborhood which
can be identified with the quotient of an open subset of a Riemannian manifold
by some finite group of isometries acting on this subset.

We omit the exact definitions for general Riemannian orbifolds, which can be
found, e.g., in \cite{Sa}, \cite{Th}, \cite{CR}, \cite{We}, because actually we will be
dealing in this article only with the special case of so-called ``good''
Riemannian orbifolds. A good Riemannian orbifold~$\O$ is the quotient of
a Riemannian manifold $(M,g)$ by some group of isometries $\Gamma$
which acts effectively and properly discontinuously on~$M$; that is, for each compact
subset $K\subset M$, the set $\{\gamma\in\Gamma\mid\gamma K\cap K\ne\emptyset\}$ is
finite. Let $p:M\to \Gamma\backslash M=\O$ be the canonical projection.
For $x\in\O$, the isotropy group $\Iso(x)$ of~$x$ is defined as the isomorphism
class of the stabilizer $\Gamma_{\tilde x}:=\{\gamma\in\Gamma\mid \gamma\tilde x=\tilde x\}$
of~$\tilde x$ in~$\Gamma$,
where $\tilde x$ is any point in the preimage $p\inv(x)\subset M$ of~$x$.
Note that $\Iso(x)$ is well-defined because for any $\tilde x'\in p\inv(x)$
the groups $\Gamma_{\tilde x}$ and $\Gamma_{\tilde x'}$ are conjugate in~$\Gamma$.
By abuse of notation we will sometimes call $\Gamma_{\tilde x}$ (instead of
its isomorphism class) the isotropy group of $x=p(\tilde x)$.
If $\Iso(x)$ is nontrivial then $x$ is called a singular point of~$\O$,
and the (finite) number $\#\Iso(x)$ is called its isotropy order.

The space $C^\infty(\O)$ of smooth functions on a good Riemannian orbifold
$\O=\Gamma\backslash M$ may be defined as the space $C^\infty(M)^\Gamma$ of
$\Gamma$-invariant smooth functions on~$M$. Similarly, smooth $k$-forms on~$\O$
are defined as $\Gamma$-invariant smooth $k$-forms on~$M$. Since the Laplace operator
$\Delta_g$ on $(M,g)$ commutes with isometries and thus preserves $\Gamma$-invariance,
it preserves the space $C^\infty(\O)$, and its restriction to this space is called
the Laplace operator on functions on~$\O$.
Similarly, the Laplace operator on $k$-forms
on~$\O$ is the restriction of $dd^*+d^*d:\Omega_k(M)\to\Omega_k(M)$ to the space of
$\Gamma$-invariant $k$-forms.
Again, these notions can be suitably defined also on general Riemannian orbifolds
and coincide with the given ones on good Riemannian orbifolds. On every compact connected
Riemannian orbifold the Laplace operator on functions has a discrete spectrum
of eigenvalues $0=\lambda_0<\lambda_1\le\lambda_2\le\ldots\to\infty$ with finite
multiplicities; see~\cite{Ch}. For a good orbifold $\O$ as above,
the eigenspace $E_\lambda(\O)$ associated with the eigenvalue~$\lambda$ of~$\O$
is canonically identified with the subspace of $\Gamma$-invariant elements
of the space $E_\lambda(M,g)$ of eigenfunctions associated with this eigenvalue on $(M,g)$.
Two compact Riemannian orbifolds are called isospectral
if they have the same spectrum. 

To which extent does the Laplace spectrum determine the geometry of a compact
Riemannian orbifold, and,
in particular, the structure of its singularities? There exist some positive
results in this direction. An important general observation is that a compact
Riemannian orbifold which is not a manifold (i.e., has singular points) can
never be isospectral to a Riemannian manifold with which it shares a common
Riemannian covering. This is shown in~\cite{GR} using an asymptotic
expansion by H.~Donnelly of the heat trace for good compact Riemannian
orbifolds; his result was made more explicit and generalized to
non-good orbifolds in~\cite{DGGW}.
It is not known whether the statement concerning nonisospectrality
of manifolds and orbifolds remains
true without the condition of a common Riemannian covering.
E.~Dryden and A.~Strohmaier showed
that on oriented compact hyperbolic orbifolds
in dimension two, the spectrum completely determines the types and numbers of
singular points \cite{DS}.
Independently, this had also been shown by the first author together with
P.G.~Doyle (unpublished).
By a result of E.~Stanhope, only finitely many isotropy groups can occur
in a family of isospectral orbifolds satisfying a uniform lower bound on the
Ricci curvature \cite{St}.
On the other hand, N.~Shams, E.~Stanhope, and D.~Webb have constructed arbitrarily
large (finite) families of mutually isospectral Riemannian orbifolds such that
each of these contains an isotropy group which does not occur in any of the other
orbifolds of the family \cite{SSW}. More precisely, for the maximal isotropy
orders occurring in the orbifolds of such a family, the corresponding isotropy
groups all have the same order, but are mutually nonisomorphic. A natural
question arising in this context is whether it might be possible that
two isospectral orbifolds have maximal isotropy groups which are not only nonisomorphic
but even of different size. The only previously known examples of this kind
concerned pairs of orbifolds with disconnected topology \cite{DR:2006}.
The present paper, however, exhibits several kinds of examples of isospectral
connected orbifolds with different maximal isotropy orders; thus, using a
popular formulation: You cannot hear the maximal isotropy order of an orbifold.

The paper is organized as follows:

In Section~\ref{sec:sunada} we recall B\'erard's, Ikeda's and Pesce's
versions of the Sunada theorem
and apply it to obtain a general construction of pairs of isospectral
orbifolds with different maximal isotropy orders (Theorem~\ref{thm:genmax},
Corollary~\ref{cor:biinv}), as well as some
explicit examples. In this approach, the orbifolds arise as quotients of
Riemannian manifolds which are locally isometric to a compact Lie group
with a biinvariant metric, or, more generally, to a homogeneous space.

In Section~\ref{sec:flat} we recall some formulas developed by R.~Miatello and the first author concerning the spectrum of flat manifolds and orbifolds.
We use these to obtain several isospectral pairs of compact
flat $3$-dimensional orbifolds, among these also a pair with different
maximal isotropy orders (Example~\ref{ex:flat1}).
In another example (Example~\ref{ex:flat2}), the maximal isotropy
groups are of the same size but not isomorphic, as in the examples by
Shams, Stanhope, and Webb~\cite{SSW}. Moreover, the sets of singular points
of maximal isotropy order have different dimension in the two orbifolds.
Example~\ref{ex:flat5} is another example of this kind.
In Examples~\ref{ex:flat3} and \ref{ex:flat4}, all nontrivial isotropy
groups are isomorphic, but again the topology of the singular sets is different.
These two examples are obtained by the classical Sunada construction. Their
existence within the context of flat $3$-dimensional orbifolds is interesting
because it is known \cite{RC}, \cite{DR:2004} that there do not exist nontrivial pairs
of Sunada isospectral flat manifolds in dimension three. See~\cite{We}
for a more detailed treatment of some of the examples in this section.

The first author would like to thank
Humboldt-Universit\"at zu Berlin, and especially Dorothee Schueth, for the great hospitality
during his one year stay there.

\Section{Sunada isospectral orbifolds}
\label{sec:sunada}

\noindent
The famous Sunada theorem \cite{Su} gives a general method for constructing
isospectral manifolds and orbifolds. In order to formulate it and the versions
given by P.~B\'erard and A.~Ikeda which we will use here, one needs the
notion of almost conjugate subgroups.

\begin{definition}
\label{def:almost}\
Let $G$ be a group. Two finite subgroups $\Gamma_1$, $\Gamma_2$ of~$G$ are
called \emph{almost conjugate in~$G$} if each conjugacy class
$[b]_G$ in $G$ intersects $\Gamma_1$ and $\Gamma_2$ in the same number
of elements: $\#([b]_G\cap\Gamma_1)=\#([b]_G\cap\Gamma_2)$.
\end{definition}

The classical version of the Sunada theorem says that if $G$ is a finite
group acting by isometries on a compact Riemannian manifold $(M,g)$, and if
$\Gamma_1$ and $\Gamma_2$ are almost conjugate subgroups of~$G$ acting without
fixed points on~$M$, then the quotient
manifolds $\Gamma_1\backslash M$, $\Gamma_2\backslash M$, each endowed
with the metric induced by~$g$, are isospectral. If one drops the condition
that $\Gamma_1$ and $\Gamma_2$ act without fixed points then the statement
remains true in the context of Riemannian orbifolds, as shown by P.~B\'erard~\cite{Be}.
Finally, A.~Ikeda \cite{Ik} showed
that the Sunada theorem still holds in the case that
$G$ is the whole (necessarily compact) group of isometries of $(M,g)$,
or any subgroup of the latter (as his proof allows).
Although he did not formulate this result for orbifolds, the proof he 
gives in the manifold context carries over verbatim to the orbifold case.
Independently, H.~Pesce \cite{Pe} had already given a version of the Sunada
theorem for compact, not necessarily finite~$G$,
with a slightly different (but equivalent)
formulation of the almost conjugacy condition in representation theoretic
terms. Thus, one has the following theorem (which can also be interpreted
as a special case of a much more general result by D.~DeTurck and C.~Gordon
\cite{DG}):

\begin{theorem}[\cite{Su}, \cite{Be}, \cite{Ik}, \cite{Pe}]
\label{thm:sunada}\
Let $(M,g)$ be a compact Riemannian manifold, and let $G$ be a
group which acts by isometries on $(M,g)$. If $\Gamma_1$ and $\Gamma_2$ are
two finite subgroups which are almost conjugate in~$G$, then the compact
Riemannian orbifolds $\O_1:=\Gamma_1\backslash M$ and $\O_2:=\Gamma_2\backslash M$
are isospectral.
\end{theorem}

Note that we have not assumed effectiveness of the action of the $\Gamma_i$ on~$M$.
However, by identifying $\Gamma_i\backslash M$ with $\tilde\Gamma_i\backslash M$,
where $\tilde\Gamma_i$ is the quotient of~$\Gamma_i$ by the kernel of its action,
this orbifold is again seen to be a good Riemannian orbifold in the sense of the
introduction.

We briefly sketch Ikeda's particularly simple proof of Theorem~\ref{thm:sunada}:
Since $G$ acts by isometries,
its canonical action on $C^\infty(M)$ commutes with the Laplace operator $\Delta_g$;
in particular, it
preserves the corresponding eigenspaces $E_\lambda(M,g)$. Fix $\lambda$, let
$V:=E_\lambda(M,g)$, and denote the action of $G$ on $V$ by $\rho$. Note that
$V$ is finite dimensional since $M$ is compact.
We have to show that for $i=1,2$, the $\Gamma_i$-invariant
subspaces $V^{\Gamma_i}\cong E_\lambda(\O_i)$ of~$V$ have the same dimension.
But this dimension is the trace of the projection operator $(\#\Gamma_i)\inv
\sum_{\gamma\in\Gamma_i}\rho_\gamma$; it is thus equal to
$(\#\Gamma_i)\inv\sum_{\gamma\in\Gamma_i}\tr(\rho_\gamma)$. Since there
exists a bijection from $\Gamma_1$ to $\Gamma_2$ which preserves conjugacy
classes in $G$, and thus traces, the two numbers are indeed the same for $i=1,2$.

\begin{remark}
\label{rem:sunadaforms}\
Sunada-isospectral orbifolds (i.e.,
isospectral orbifolds arising from Theorem~\ref{thm:sunada}) are actually
isospectral on $k$-forms for all~$k$; see the articles cited above. In fact,
the above proof goes through without change if one replaces
smooth functions by smooth $k$-forms.
\end{remark}

If $\Gamma_1$ and $\Gamma_2$ are not only almost conjugate, but conjugate
in~$G$ then the situation becomes trivial; in fact, if $\Gamma_2 = a\Gamma_1a\inv$
for some $a\in G$ then $a:M\to M$ induces an isometry between the Riemannian orbifolds
$\Gamma_1\backslash M$ and $\Gamma_2\backslash M$. Fortunately there exist many
triples $(G,\Gamma_1,\Gamma_2)$ where the $\Gamma_i$ are almost
conjugate, but not conjugate in~$G$. One example which we are going to use
is the following:

\begin{example}
\label{ex:so6}\
Let $G:=\SO(6)$. Writing diagonal matrices in~$G$ as the vectors of their
entries on the diagonal, define
\begin{align*}
\Gamma_1:=\{&(1,1,1,1,1,1), (-1,-1,-1,-1,-1,-1),\\
&(-1,-1,1,1,1,1), (-1,1,-1,1,1,1), (1,-1,-1,1,1,1),\\
&(-1,1,1,-1,-1,-1), (1,-1,1,-1,-1,-1), (1,1,-1,-1,-1,-1)\},\\
\Gamma_2:=\{&(1,1,1,1,1,1), (-1,-1,-1,-1,-1,-1),\\
&(-1,-1,1,1,1,1), (1,1,-1,-1,1,1), (1,1,1,1,-1,-1),\\
&(-1,-1,-1,-1,1,1), (-1,-1,1,1,-1,-1), (1,1,-1,-1,-1,-1)\}.
\end{align*}
Obviously there is a bijection from $\Gamma_1$ to $\Gamma_2$ preserving conjugacy
classes in $G$; thus the two subgroups are almost conjugate in $G$. (Actually,
the two groups can be seen to be almost conjugate by elements of the group~$A_6$
of even permutation matrices in $G$, and thus almost conjugate \emph{in} the
finite subgroup of $G$ generated by $\Gamma_1\cup\Gamma_2\cup A_6$.)
This example corresponds to a certain pair of linear
codes in $\Z_2^6$ with the same weight enumerator, mentioned in~\cite{CS}.
The groups $\Gamma_1$ and $\Gamma_2$ are not conjugate in $G=\SO(6)$ because
$\Gamma_1$ has a four-element subgroup acting as the identity on some
three-dimensional
subspace of~$\R^3$ (namely,  on $\spann\{e_4,e_5,e_6\}$), while no four-element
subgroup of~$\Gamma_2$ acts as the identity on any three-dimensional subspace
of~$\R^3$.

\end{example}

The following observation is the main point of this section:

\begin{theorem}
\label{thm:genmax}\
Let $G$ be a compact Lie group and $H$ be a closed Lie subgroup of~$G$.
Choose a left invariant Riemannian metric on~$G$ which is also right invariant
under~$H$. Let $g$ denote the corresponding
Riemannian metric on the quotient manifold $M:=G/H$ such that the canonical
projection $G\to G/H$ becomes a Riemannian submersion. Let $\Gamma_1$ and $\Gamma_2$
be two finite subgroups of~$G$ which are almost conjugate in~$G$.
\begin{itemize}
\item[(i)] 
The compact Riemannian orbifold quotients $\O_1:=\Gamma_1\backslash M$ and
$\O_2:=\Gamma_2\backslash M$ of $(M,g)$ are isospectral.
\item[(ii)] Let $m(\Gamma_i,H):=\max_{a\in G}\#(\Gamma_i\cap aHa\inv)$
and $n(\Gamma_i,H):=\#(\Gamma_i\cap\bigcap_{a\in G} aHa\inv)$ for $i=1,2$.
Then $m(\Gamma_i,H):n(\Gamma_i,H)$ is the maximal isotropy order of singular points in~$\O_i$.
Moreover, $n(\Gamma_1,H)=n(\Gamma_2,H)$. In particular, if $m(\Gamma_1,H)\ne m(\Gamma_2,H)$
then $\O_1$ and $\O_2$ have different maximal isotropy orders.
\end{itemize}
\end{theorem}

\begin{proof}
(i) This follows from Theorem~\ref{thm:sunada} because $G$ acts by isometries
on the homogeneous space $(M,g)=(G/H,g)$.

(ii) Let $a\in G$. Then the stabilizer in $\Gamma_i$ of the point $aH\in M$
is the group $\{\gamma\in\Gamma_i\mid \gamma aH = aH\} = \{\gamma\in\Gamma_i\mid
\gamma\in aHa\inv\}$; that is,
\begin{equation}
\label{isogr}
(\Gamma_i)_{aH} = \Gamma_i\cap aHa\inv.
\end{equation}
Moreover, the kernel of the action of~$\Gamma_i$ on $G/H$ is $\Gamma_i\cap\bigcap_{a\in G}
aHa\inv$.
This implies the formula for the maximal isotropy orders. For the statement about
the numbers $n(\Gamma_i,H)$ let $\Phi:\Gamma_1\to\Gamma_2$ be a bijection which
preserves $G$-conjugacy classes. Note that $N:=\bigcap_{a\in G}aHa\inv$ is a
normal subgroup of~$G$. Hence $\Phi$ restricts to a bijection from $\Gamma_1\cap N$
to $\Gamma_2\cap N$.
\end{proof}

\begin{corollary}
\label{cor:biinv}\
Let $G$ be a compact Lie group and $\Gamma_1$, $\Gamma_2$ be two almost conjugate,
non-conjugate finite subgroups of~$G$. Choose a biinvariant metric on~$G$, and denote
the induced metric on the quotient manifold $M:=G/\Gamma_1$ by~$g$.
Then the compact Riemannian orbifold quotients $\O_1:=\Gamma_1\backslash M$ and $\O_2:=
\Gamma_2\backslash M$ of $(M,g)$ are isospectral and have different maximal isotropy
orders.
\end{corollary}

\begin{proof}
This follows immediately from Theorem~\ref{thm:genmax} with $H:=\Gamma_1$. In fact,
we have $m(\Gamma_1,\Gamma_1)=\#\Gamma_1=\#\Gamma_2$; if this were equal
to $m(\Gamma_2,\Gamma_1)$ then $\Gamma_1$ and $\Gamma_2$ would be conjugate
by some $a\in G$, contradicting the hypothesis.
\end{proof}

\begin{example}
\label{ex:so6a}\
The following is an example for Theorem~\ref{thm:genmax} not arising from the corollary.
Let $G,\Gamma_1,\Gamma_2$ be the groups from Example~\ref{ex:so6}.
Let $H\cong\SO(3)$ be the subgroup of $G$ consisting of matrices
of the form
$$
\left(\begin{matrix} A&0\\0&I_3 \end{matrix}\right),
$$
where $I_3$ denotes the unit element in $\SO(3)$. Then $M:=G/H=\SO(6)/\SO(3)$
is the Stiefel manifold $V_{6,3}$ of orthonormal $3$-frames in euclidean $\R^6$; the
point $aH\in M$ corresponds to the $3$-frame formed by the three last
column vectors of the matrix $a\in\SO(6)$.
Note that $G$ acts effectively on~$M$.
Choose a biinvariant metric on $\SO(6)$ (or any left invariant metric which
is also right invariant under~$H$) and endow $M$ with the corresponding
homogeneous metric. By Theorem~\ref{thm:genmax}, the compact Riemannian
orbifolds $\O_1:=\Gamma_1\backslash M$ and $\O_2:=\Gamma_2\backslash M$ are isospectral.
Moreover, the point $eH\in M$ (corresponding to the orthonormal $3$-frame
$(e_4,e_5,e_6)$ in $V_{6,3}$, where $e_i$ denotes the $i$-th standard unit vector)
is stabilized by four elements in $\Gamma_1$, namely, the elements of $\Gamma_1\cap H$
(recall~(\ref{isogr})).
The same point is also stabilized by some two-element subgroup
of $\Gamma_2$. On the other hand, no four-element subgroup of $\Gamma_2$
stabilizes any point in~$M$: Such a point would have to correspond to an orthonormal
$3$-frame each of whose vectors is contained in the intersection of the $1$-eigenspaces
of the group elements; but for each four-element subgroup of $\Gamma_1$ this intersection
is at most two-dimensional. Since obviously no point in $M$ (not even any single unit vector
in~$\R^6$) is stabilized by the whole group $\Gamma_1$, we see
that $\O_1$ has maximal isotropy order four, while $\O_2$ has maximal isotropy order two.
In the notation of Theorem~\ref{thm:genmax}, $m(\Gamma_1,H)=4$, $m(\Gamma_2,H)=2$,
and $n(\Gamma_1,H)=n(\Gamma_2,H)=1$.
\end{example}

\begin{example}
\label{ex:so6b}\
Let $G,\Gamma_1,\Gamma_2$ again be as in the previous example, and let $g$ be the
Riemannian metric on $M:=G/\Gamma_1$ induced by a biinvariant metric on $G=\SO(6)$.
Then the Riemannian orbifold quotients $\O_1:=\Gamma_1\backslash M$ and
$\O_2:=\Gamma_2\backslash M$
of $(M,g)$ are isospectral and have different maximal isotropy orders by
Corollary~\ref{cor:biinv}.

More precisely, the maximal isotropy order of singular points in $\O_1$
is $m(\Gamma_1,\Gamma_1):2=4$, while in $\O_2$ it is $m(\Gamma_2,\Gamma_1):2=2$.
In fact, $N:=\bigcap_{a\in G}a\Gamma_1a\inv\subset\Gamma_1\cap\Gamma_2$
is the sub\-group~$\{\pm I_6\}$ of order~$2$, and we have $m(\Gamma_2,\Gamma_1)=4$
because a four-element subgroup of~$\Gamma_1$ which contains
$-I_6$ is conjugate by some $a\in G$ (for example, a permutation matrix)
to a subgroup of~$\Gamma_2$.
\end{example}

\begin{example}
\label{ex:so6c}\
Another variation of the above examples, but \emph{not} leading
to different maximal isotropy orders, is obtained by letting $G$
act canonically on the standard unit sphere $(M,g):=S^5$;
in our above approach, this corresponds to letting $H:=\SO(5)$.
As one immediately sees, the isotropy group of maximal order in
$\O_i:=\Gamma_i\backslash S^5$ is isomorphic to $\Z_2\times\Z_2$ for both $i=1,2$.
Nevertheless it is possible to distinguish
between $\O_1$ and $\O_2$ by using the \emph{topology} of the
set $\S_i\subset\O_i$ of singularities with maximal isotropy orders,
that is, the image in $\O_i$ of the set of points in $S^5$ whose
stabilizer in $\Gamma_i$
consists of four elements: The set $\S_1$ is the disjoint
union of one copy of $\R P^2$ (the image of the
unit sphere in $\spann\{e_4,e_5,e_6\}$) and of three points
(the images of $\pm e_1$, $\pm e_2$, and $\pm e_3$). The set $\S_2$,
in contrast, is the disjoint union of
three copies of $S^1$ (the images of the unit spheres
in $\spann\{e_1,e_2\}$, $\spann\{e_3,e_4\}$,
and $\spann\{e_5,e_6\}$).
\end{example}

\begin{remark}
\label{rem:topsph}\
(i) The fact that the topological structure of certain singular strata can be different
in isospectral orbifolds has also been shown in~\cite{SSW}; a new feature in Example~\ref{ex:so6c}
is that this concerns the set of points of \emph{maximal} isotropy order.
We will reencounter the analogous situation in certain isospectral
pairs of flat $3$-dimensional orbifolds; see Examples
\ref{ex:flat2}, \ref{ex:flat3}, \ref{ex:flat4}, and \ref{ex:flat5}.

(ii) It is easy to see that for almost conjugate pairs $\Gamma_1,\Gamma_2$ of
\emph{diagonal} subgroups of $\SO(n)$, necessarily containing only $\pm1$ as entries (as the pair used in the
above examples), the corresponding actions
on $S^{n-1}\cong\SO(n)/\SO(n-1)$ will always have the same maximal isotropy order
(and isomorphic maximal isotropy groups $\Z_2^k$ for some $k$).
We do not know whether there exist pairs of almost conjugate finite subgroups
$\Gamma_1$ and $\Gamma_2$ of $\SO(n)$ which satisfy $m(\Gamma_1,\SO(n-1))\ne
m(\Gamma_2,\SO(n-1))$ and
would thus yield, by Theorem~\ref{thm:genmax},
isospectral \emph{spherical} orbifolds with different maximal isotropy orders.
\end{remark}

\begin{remark}
Once one has a pair of isospectral compact Riemannian
orbifolds $\O_1$, $\O_2$ with different maximal isotropy orders, then
one immediately obtains for each $m\in\N$ a family of $m+1$ mutually isospectral
Riemannian orbifolds $\V_0,\ldots,\V_m$ with pairwise different maximal
isotropy orders; one just defines $\V_i$ as the Riemannian product
of $i$ times $\O_1$ and $m-i$ times $\O_2$.
The Riemannian product of two \emph{good} Riemannian orbifolds (as are
all orbifolds in our examples) $\O=\Gamma\backslash M$ and $\O'=\Gamma'
\backslash M'$ of $(M,g)$, resp.~$(M',g')$, is defined as $(\Gamma\times
\Gamma')\backslash (M\times M')$, where $M\times M'$ is endowed with
the Riemannian product metric associated with $g$ and~$g'$.
\end{remark}

\Section{Isospectral flat orbifolds in dimension three}
\label{sec:flat}

\noindent
A Riemannian orbifold~$\O$ is called flat if each point in~$\O$ has a
neighborhood which is the quotient of an open subset of~$\R^n$, endowed
with the euclidean metric, by a finite group of Riemannian isometries.
It can be shown that every flat orbifold is good \cite{Th}; hence,
it is the quotient of a flat Riemannian manifold by some group of isometries
acting properly discontinuously.

Let us recall some facts from the theory of quotients of standard
euclidean space~$(\R^n,g)$ by groups of isometries; see \cite{Wo}.
The isometry group $I(\R^n,g)$ is the semidirect product $O(n)\ltimes\R^n$
consisting of all transformations $BL_b$ with $B\in O(n)$ and $b\in\R^n$,
where $L_b$ is the translation $x\mapsto x+b$ of $\R^n$.
Note that
\begin{equation}
\label{irnmult}
L_bB=BL_{B\inv b},\mbox{ }BL_bB\inv=L_{Bb},\mbox{ and }(BL_b)\inv=B\inv L_{-Bb}.
\end{equation}
The compact-open topology on $I(\R^n,g)$ coincides with the canonical
product topology on $O(n)\times\R^n$. A subgroup $\Gamma$ of $I(\R^n,g)$
acts properly discontinuously with compact quotient on~$\R^n$ if and only if
it is discrete and cocompact in $I(\R^n,g)$. Such a group is called a
crystallographic group. If, in addition, $\Gamma$ is torsion-free,
then it acts without fixed points on~$\R^n$, and $\Gamma\backslash\R^n$ is a flat
Riemannian manifold. Conversely, every compact flat Riemannian manifold is isometric
to such a quotient. If the condition that $\Gamma$ be torsion-free is dropped
then $\Gamma\backslash\R^n$ is a compact good Riemannian orbifold which is flat.
Conversely, if $\O$ is any compact flat Riemannian orbifold (and is thus,
as mentioned above, a good orbifold), then there exists a crystallographic
group $\Gamma\subset I(\R^n,g)$ such that $\O$ is isometric to $\Gamma\backslash\R^n$.

If $\Gamma$ is a crystallographic group acting on $\R^n$ then the translations
in $\Gamma$ form a normal, maximal abelian subgroup~$L_\Lambda$ where $\Lambda$
is a cocompact lattice in~$\R^n$; the quotient group $\bar\Gamma:=\Gamma/L_\Lambda$
is finite.
The flat torus $T_\Lambda:=L_\Lambda\backslash\R^n$
covers $\O:=\Gamma\backslash\R^n$ because
$L_\Lambda$ is normal in~$\Gamma$. More precisely, we have
$\O\cong\bar\Gamma\backslash T_\Lambda$,
where $\gamma L_\Lambda\in\bar\Gamma$ acts on $T_\Lambda$ as the map
$\bar\gamma:T_\Lambda\to T_\Lambda$ induced by $\gamma:\R^n\to\R^n$.
Let $F\subset O(n)$ be the image of the canonical
projection from $\Gamma\subset O(n)\ltimes\R^n$ to~$O(n)$. This projection
has kernel~$L_\Lambda$; thus we have $F\cong\bar\Gamma$.

Let $k\in\{0,\ldots,n\}$.
For $\mu\ge0$ let $H_{k,\mu}(T_\Lambda)$ denote the space of smooth $k$-forms on $T_\Lambda$ which are eigenforms
associated with the eigenvalue $4\pi^2\mu$. Then the multiplicity of
$4\pi^2\mu$ as an eigenvalue for the Laplace operator on $k$-forms on
the Riemannian orbifold $\O=\Gamma\backslash\R^n=\bar\Gamma\backslash T_\Lambda$ equals the
dimension of the subspace
$$
H_{k,\mu}(T_\Lambda)^{\bar\Gamma} = \{\omega\in H_{k,\mu}(T_\Lambda)\mid
\bar\gamma^*\omega=\omega\;
\forall\bar\gamma\in\bar\Gamma\}
$$
(which might be zero). This dimension can be computed using the formula
from the following theorem.

\begin{theorem}[\cite{MR:2001}, \cite{MR:2002}]
\label{thm:dimformula}\
Let $d_{k,\mu}(\Gamma):=\dimm H_{k,\mu}(T_\Lambda)^{\bar\Gamma}$. Then
$$
d_{k,\mu}(\Gamma)=(\#F)\inv\sum_{B\in F}\tr_k(B)e_{\mu,B}(\Gamma),\mbox{\ \ where \;}
e_{\mu,B}(\Gamma):=\sum_{{v\in\Lambda^*,\|v\|^2=\mu}\atop{Bv=v}}e^{2\pi i\<v,b\>}
$$
with $b$ chosen such that $BL_b\in\Gamma$,
the trace of $B$ acting on the $\binom nk$-dimensional
space of alternating $k$-linear forms on~$\R^n$ as pullback by~$B\inv$
is denoted by $\tr_k(B)$, and where $\Lambda^*:=\{v\in\R^n\mid \<v,\lambda\>\in\Z\;
\forall\lambda\in\Lambda\}$ is the dual lattice associated with~$\Lambda$.
\end{theorem}

\begin{notrems}\ \\
(i) Note that $\tr_0(B)=1$ and $\tr_1(B)=\tr(B\inv)=\tr(\trans B)=\tr(B)$
for all $B\in O(n)$.\\
(ii) For $k=0$ we write $d_\mu:=d_{0,\mu}$. Thus $d_\mu(\Gamma)$ will be
the multiplicity of $4\pi^2\mu$ as an eigenvalue for the Laplace operator
on functions on~$\O$.
\end{notrems}

The following is an example of two isospectral flat three-dimensional orbifolds with
different maximal isotropy orders.

\begin{example}
\label{ex:flat1}\
Let $\Lambda$ be the lattice $2\Z\times2\Z\times\Z$ in~$\R^3$. Define
$$
\tau:=\left(\begin{matrix} 0&-1&0\\1&0&0\\0&0&1\end{matrix}\right),\mbox{ }
\chi_1:=\left(\begin{matrix} 1&0&0\\0&-1&0\\0&0&-1\end{matrix}\right),\mbox{ }
\chi_2:=\left(\begin{matrix} -1&0&0\\0&1&0\\0&0&-1\end{matrix}\right),\mbox{ }
\chi_3:=\left(\begin{matrix} -1&0&0\\0&-1&0\\0&0&1\end{matrix}\right)
$$
and
$$
b_1:=e_1,\mbox{ }b_2:=0,\mbox{ }b_3:=-e_1\in\R^3.
$$
Let $\Gamma_1$ be the subgroup of $I(\R^3)$ generated by $L_\Lambda$ and $\tau$, and let
$\Gamma_2$ be generated by $L_\Lambda$ and the maps $\rho_j:=\chi_j\circ L_{b_j}$
($j=1,2,3$). Using~(\ref{irnmult}) one easily checks that
$$
\Gamma_1=\{\tau^j L_\lambda\mid j\in\{0,1,2,3\}, \lambda\in\Lambda\}\mbox{ and }
\Gamma_2=\{\rho_j L_\lambda\mid j\in\{0,1,2,3\}, \lambda\in\Lambda\},
$$
where $\rho_0:=\Id$.
Since these are discrete and cocompact subgroups of $O(3)\ltimes\R^3$, we obtain two compact
flat orbifolds
$$
\O_1:=\Gamma_1\backslash\R^3,\mbox{ and }\O_2:=\Gamma_2\backslash\R^3.
$$
It is not difficult to see that the unit cube $[0,1]^3$ is a fundamental domain
for the action of $\Gamma_1$, resp.~$\Gamma_2$, on~$\R^3$, and that the identifications
on the sides are as given in the following two figures, where the top and bottom
sides are identified by the vertical translation~$L_{e_3}$.

\begin{figure}[!htb]
\centerline{\mbox{\includegraphics*[scale=.4]{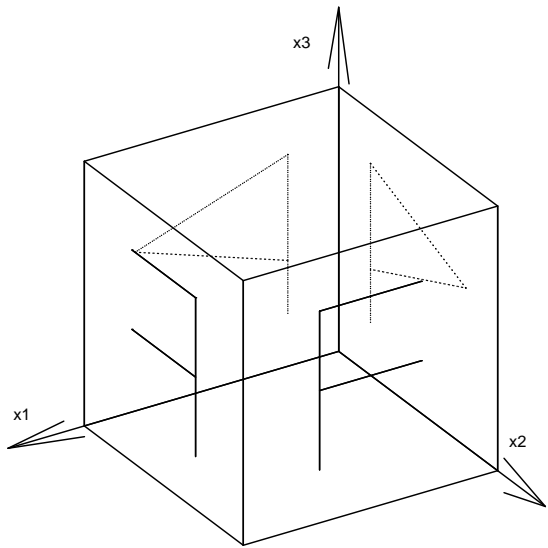}}}
\caption{The underlying space of $\O_1$ as a quotient of the unit cube}
\label{fig:cube1}
\end{figure}

\begin{figure}[htb]
\centerline{\mbox{\includegraphics*[scale=.4]{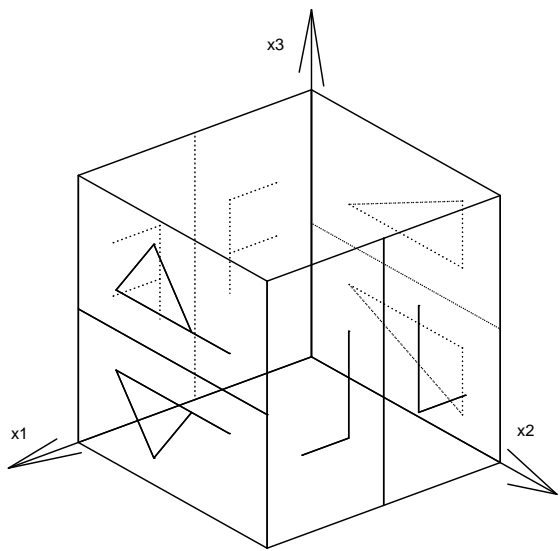}}}
\caption{The underlying space of $\O_2$ as a quotient of the unit cube}
\label{fig:cube2}
\end{figure}

In Figure~\ref{fig:cube1}, describing~$\O_1$,
the element $\tau\in\Gamma_1$ accounts for the side identification denoted
by~$\mathsf{P}$, and $\tau L_{-2e_2}\in\Gamma_1$ for the one denoted by~$\mathsf{F}$.
Note that $\O_1$ is actually the Riemannian product of a two-dimensional
so-called $442$-orbifold and a circle of length one. (A $442$-orbifold has
two cone points of order~$4$ and one cone point of order~$2$.)
In Figure~\ref{fig:cube2} which describes~$\O_2$, the elements of~$\Gamma_2$ which account
for the side identifications denoted by $\mathsf F$, $\mathsf L$, $\mathsf{\Delta}$,
$\mathsf{P}$ are $\rho_3$, $\rho_3 L_{-2e_2}$, $\rho_2 L_{-e_3}$, and $\rho_2 L_{-2e_1-e_3}$,
respectively.

\textit{Isotropy groups:}
It is clear that the isotropy groups both in $\O_1$ and $\O_2$ can have \emph{at most}
order four because $L_\Lambda$ has index four in~$\Gamma_i$ and a point in~$\R^3$
cannot be fixed simultaneously by two isometries that differ by a nontrivial translation.

Since $\tau$ is a quarter rotation around the axis spanned by~$e_3$, the
four-element subgroup $\{\Id,\tau,\tau^2,\tau^3\}\cong\Z_4$ of~$\Gamma_1$
pointwise fixes the
edge $\{(0,0,x_3)\mid 0\le x_3\le1\}$ of the fundamental cube; thus $\O_1$
has maximal isotropy order four. The other points in the fundamental domain
with nontrivial stabilizer in~$\Gamma_1$ are $\{(1,1,x_3)\mid 0\le x_3\le1\}$,
pointwise fixed by the four-element group generated by $\tau L_{-2e_2}$,
and $\{(0,1,x_3)\mid0\le x_3\le1\}$ (identified with $\{(1,0,x_3)\mid0\le x_3\le1\}$
via the identifications marked $\mathsf{F}$ or~$\mathsf{P}$ in Figure~\ref{fig:cube1}),
pointwise fixed by $\{\Id, \tau^2L_{-2e_2}\}\cong\Z_2$.
So the singular set in~$\O_1$ consists of three copies of~$S^1$, each of length one,
two of them with isotropy group~$\Z_4$ and one with~$\Z_2$. (Of course, these three
components correspond to the three cone points of the $442$-orbifold mentioned above.)

In $\O_2$ there are no points with isotropy order four. Otherwise, there would
have to exist a point in~$\R^3$ fixed by three elements of the form $\rho_1L_\lambda$,
$\rho_2L_\mu$, $\rho_3L_\nu$ with $\lambda,\mu,\nu\in\Lambda$. But $(\rho_1L_\lambda)^2
=(\chi_1L_{e_1+\lambda})^2=L_{(\chi_1\inv+\Id)(e_1+\lambda)}$. In order to fix a point,
this translation would have to be trivial; in particular, the first coordinate
of $e_1+\lambda$ would have to vanish. This contradicts $\lambda\in\Lambda$.
Thus, the points in $\O_2$ which do have nontrivial isotropy all have isotropy group~$\Z_2$.
The singular set in~$\O_2$ consists of four copies of~$S^1$: Two of length two, corresponding
to the horizontal edges and middle segments in the faces of the fundamental cube marked by $\mathsf{P}$ and~$\mathsf{\Delta}$ in Figure~\ref{fig:cube2},
and two of length one, corresponding
to the middle vertical segments on the faces marked by $\mathsf{L}$ and~$\mathsf{F}$.

\textit{Isospectrality:}
Let $\mu\ge0$. The space of eigenfunctions associated with the eigenvalue $4\pi^2\mu$
on $\O_i$ has dimension $d_\mu(\Gamma_i)$ ($i=1,2$) which we compute using
Theorem~\ref{thm:dimformula} with $k=0$. We have $F_1=\{\Id,\tau,\tau^2,\tau^3\}$ and $F_2=
\{\Id,\chi_1,\chi_2,\chi_3\}$. Obviously, $d_0(\Gamma_i)=1$ for both $i=1,2$. Let $\mu>0$.
For $B=\Id$, we get $e_{\mu,\Id}(\Gamma_i)=
\#\{v\in\Lambda^*\mid\|v\|^2=\mu\}=:e_{\mu,\Id}$ for both $i=1,2$.
Note that $\Lambda^*=\frac12\Z\times\frac12\Z\times\Z$.
The only vectors of length~$\sqrt\mu$ in~$\R^3$ which are fixed by some nontrivial
element of $F_i$ are $\pm\sqrt\mu e_3$ for $i=1$ and $\pm\sqrt\mu e_j$ ($j=1,2,3$)
for $i=2$.
Therefore, if $\sqrt\mu\notin\frac12\N$ then no $v\in\Lambda^*$ of
length~$\sqrt\mu$ is fixed by any nontrivial element of the $F_i$, and thus
$d_\mu(\Gamma_1)=\frac14e_{\mu,\Id}=d_\mu(\Gamma_2)$.
If $\sqrt\mu\in\N$ then
$$
e_{\mu,\tau^j}(\Gamma_1)=e^{2\pi i\<\sqrt\mu e_3,0\>}+e^{2\pi i\<-\sqrt\mu e_3,0\>}=2
$$
for $j=1,2,3$, and
\begin{align*}
e_{\mu,\chi_1}(\Gamma_2)&=e^{2\pi i\<\sqrt\mu e_1,e_1\>}+e^{2\pi i\<-\sqrt\mu e_1,e_1\>}=2,\\
e_{\mu,\chi_2}(\Gamma_2)&=e^{2\pi i\<\sqrt\mu e_2,0\>}+e^{2\pi i\<-\sqrt\mu e_2,0\>}=2,\\
e_{\mu,\chi_3}(\Gamma_2)&=e^{2\pi i\<\sqrt\mu e_3,-e_1\>}+e^{2\pi i\<-\sqrt\mu e_3,-e_1\>}=2,
\end{align*}
hence $d_\mu(\Gamma_1)=\frac14(e_{\mu,\Id}+6)=d_\mu(\Gamma_2)$.
Finally, if $\sqrt\mu\in\N_0+\frac12$ then $\pm\sqrt\mu e_3\notin\Lambda^*$ and thus
$e_{\mu,\tau^j}(\Gamma_1)=0$ for $j=1,2,3$ and $e_{\mu,\chi_3}(\Gamma_2)=0$; moreover,
\begin{align*}
e_{\mu,\chi_1}(\Gamma_2)&=e^{2\pi i\<\sqrt\mu e_1,e_1\>}+e^{2\pi i\<-\sqrt\mu e_1,e_1\>}=-2,\\
e_{\mu,\chi_2}(\Gamma_2)&=e^{2\pi i\<\sqrt\mu e_2,0\>}+e^{2\pi i\<-\sqrt\mu e_2,0\>}=2,
\end{align*}
hence $d_\mu(\Gamma_1)=\frac14e_{\mu,\Id}=d_\mu(\Gamma_2)$. We have now shown that
$d_\mu(\Gamma_1)=d_\mu(\Gamma_2)$ for every $\mu\ge0$; that is, $\O_1$ and $\O_2$ are
isospectral on functions.
\end{example}

\begin{remark}
\label{rem:flat1not1}\
The orbifolds $\O_1$ and $\O_2$ from the previous example are \emph{not}
isospectral on $1$-forms, as we
can compute by using Theorem~\ref{thm:dimformula} with $k=1$. Note that
$\tr(\Id)=3$, $\tr(\tau)=\tr(\tau^3)=1$ and $\tr(\tau^2)=\tr(\chi_j)=-1$
for $j=1,2,3$. Now consider $\mu>0$ with $\sqrt\mu\in\N$. Adjusting the
trace coefficients in the corresponding computation above, we get
$$
d_{1,\mu}(\Gamma_1)=\frac14(3e_{\mu,\Id}+2-2+2)
\ne\frac14(3e_{\mu,\Id}-2-2-2)=d_{1,\mu}(\Gamma_2).
$$
\end{remark}

In the following pair of isospectral flat orbifolds, the maximal isotropy orders
coincide, but the maximal isotropy groups are not isomorphic, similarly as in the
spherical examples from~\cite{SSW}. In contrast to those examples from~\cite{SSW},
the sets of singularities with maximal isotropy order will have different dimensions
in the two orbifolds.

\begin{example}
\label{ex:flat2}\
Let $\Lambda:=2\Z\times2\Z\times2\Z\subset\R^3$. Define
$\tau,\chi_1,\chi_2,\chi_3$ as in Example~\ref{ex:flat1}, let $\Gamma_1$ be generated
by~$L_\Lambda$ and~$\tau$, and let $\Gamma_2$ be generated by $L_\Lambda$ and
the $\rho_j:=\chi_j$ ($j=1,2,3$); note that the $\rho_j$ have no translational
parts this time.
Again we confirm, using~(\ref{irnmult}), that
$$
\Gamma_1=\{\tau^j L_\lambda\mid j\in\{0,1,2,3\}, \lambda\in\Lambda\}\mbox{ and }
\Gamma_2=\{\rho_j L_\lambda\mid j\in\{0,1,2,3\}, \lambda\in\Lambda\}
$$
(where $\rho_0:=\Id$), and we obtain two compact flat orbifolds
$\O_1:=\Gamma_1\backslash\R^3$ and $\O_2:=\Gamma_2\backslash\R^3$. This time,
$[0,1]\times[0,1]\times[0,2]$ is a fundamental domain for the action of $\Gamma_1$,
resp.~$\Gamma_2$, on~$\R^3$. The side identifications are given in Figure~\ref{fig:cuboids};
the top and bottom sides are again identified via the corresponding translation
$L_{2e_3}$.

%

\begin{figure}[!htb]
\centerline{\mbox{\includegraphics*[height=7cm]{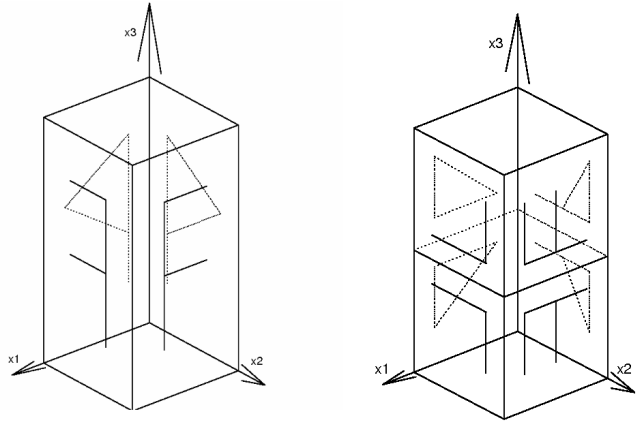}}}
\caption{The underlying spaces of $\O_1$ and $\O_2$ from Example \ref{ex:flat2}}
\label{fig:cuboids}
\end{figure}

The orbifold $\O_1$ which is pictured on the left hand side of Figure~\ref{fig:cuboids} is just a double covering of the first orbifold from the
previous example in Figure~\ref{fig:cube1}; the explanations concerning the side
identifications and
the isotropy groups are the same as before, except that now all the vertical
circles have length~$2$.
As for the right hand side of Figure~\ref{fig:cuboids}, showing $\O_2$,
the elements of~$\Gamma_2$ which account
for the side identifications denoted by $\mathsf F$, $\mathsf{\Delta}$, $\mathsf L$,
$\mathsf P$ are $\rho_1L_{-2e_2-2e_3}$, $\rho_1L_{-2e_3}$, $\rho_2 L_{-2e_1-2e_3}$,
and $\rho_2 L_{-2e_3}$, respectively.

\textit{Isotropy groups:}
One easily verifies that for $j=1,2,3$, an element $\rho_jL_\lambda\in\Gamma_2$ has fixed
points if and only if $\lambda_j=0$ (that is, $(\rho_jL_\lambda)^2=\Id$),
and in this case the fixed point set is the line
$\{-\frac12\lambda+re_j\mid r\in\R\}$. Since $\Lambda=2\Z\times2\Z\times2\Z$, the points
in~$\R^3$ with exactly two coordinates in~$\Z$ have isotropy group~$\Z_2$, while those 
in $\Z^3$ have isotropy group isomorphic to $\{\Id,\rho_1,\rho_2,\rho_3\}
\cong\Z_2\times\Z_2$. Thus (taking identifications into account),
the singular set in~$\O_2$ consists of eight
points with isotropy group $\Z_2\times\Z_2$ and of twelve
open segments of length one with isotropy group~$\Z_2$.

Since the maximal isotropy group occurring in~$\O_1$ was $\Z_4$,
the maximal isotropy orders coincide here, but the maximal
isotropy groups are nonisomorphic. Moreover, the set
of singular points with maximal isotropy has dimension one in~$\O_1$
and dimension zero in~$\O_2$.

\textit{Isospectrality:}
We continue to use the notation from the isospectrality discussion in Example~\ref{ex:flat1}
and note that now $\Lambda^*=\frac12\Z\times\frac12\Z\times\frac12\Z$.
We have $d_\mu(\Gamma_1)=\frac14e_{\mu,\Id}=d_\mu(\Gamma_2)$ if
$\sqrt\mu\notin\frac12\N$; if $\sqrt\mu\in\frac12\N$ then $e_{\mu,\tau^j}(\Gamma_1)=2$,
$e_{\mu,\rho_j}(\Gamma_2)=2$ for $j=1,2,3$, hence
$d_\mu(\Gamma_1)=\frac14(e_{\mu,\Id}+6)=d_\mu(\Gamma_2)$.
Thus $\O_1$ and $\O_2$ are isospectral on functions.
\end{example}

\begin{remark}
Similarly as in Remark~\ref{rem:flat1not1}, one shows that here
$d_{1,\mu}(\Gamma_1)\ne d_{1,\mu}(\Gamma_2)$ for $\sqrt\mu\in\frac12\N$. Thus,
$\O_1$ and $\O_2$ from Example~\ref{ex:flat2} are again \emph{not} isospectral on $1$-forms,
and, in particular, not Sunada-isospectral.
\end{remark}

The following two examples are pairs of compact flat three-dimensional orbifolds
which are Sunada-isospectral; recall that we mean by this: which arise from
Theorem~\ref{thm:sunada}. Actually, the group $G$ from the theorem
will even be finite here.
The existence of such pairs in the category of three-dimensional
flat orbifolds is noteworthy because there are no such pairs in the category of flat
three-dimensional manifolds. In fact, as shown by J.H. Conway and the first author in~\cite{RC}, there is exactly one pair, up to scaling, of isospectral flat manifolds
in dimension three. But the manifolds in that pair are not isospectral
on $1$-forms \cite{DR:2004}, and thus not Sunada-isospectral.

\begin{example}
\label{ex:flat3}\
Let $\Lambda:=\Z\times\Z\times\frac1{\sqrt2}\Z\subset\R^3$. Define
$$
\tau:=\left(\begin{matrix} 0&-1&0\\-1&0&0\\0&0&-1\end{matrix}\right)\mbox{ and }
\rho:=\left(\begin{matrix} -1&0&0\\0&-1&0\\0&0&1\end{matrix}\right).
$$
Let $\Gamma_1$ be generated
by~$L_\Lambda$ and~$\tau$, and let $\Gamma_2$ be generated by
$L_\Lambda$ and~$\rho$. Then
$$
\Gamma_1=\{\tau^j L_\lambda\mid j\in\{0,1\}, \lambda\in\Lambda\}\mbox{ and }
\Gamma_2=\{\rho^j L_\lambda\mid j\in\{0,1\}, \lambda\in\Lambda\}.
$$
Let $\O_1:=\Gamma_1\backslash\R^3$ and $\O_2:=\Gamma_2\backslash\R^3$.
A fundamental domain for
the action of $\Gamma_1$, resp.~$\Gamma_2$, on~$\R^3$, is given by
the prism of height $1/\sqrt2$ over the triangle with vertices $0$, $e_1$, $e_2$.
The side identifications are given in Figure~\ref{fig:sunada-indirect} (where once more
the top and bottom sides are identified via a vertical translation).

%

\begin{figure}[!htb]
\centerline{\mbox{\includegraphics*[height=7cm]{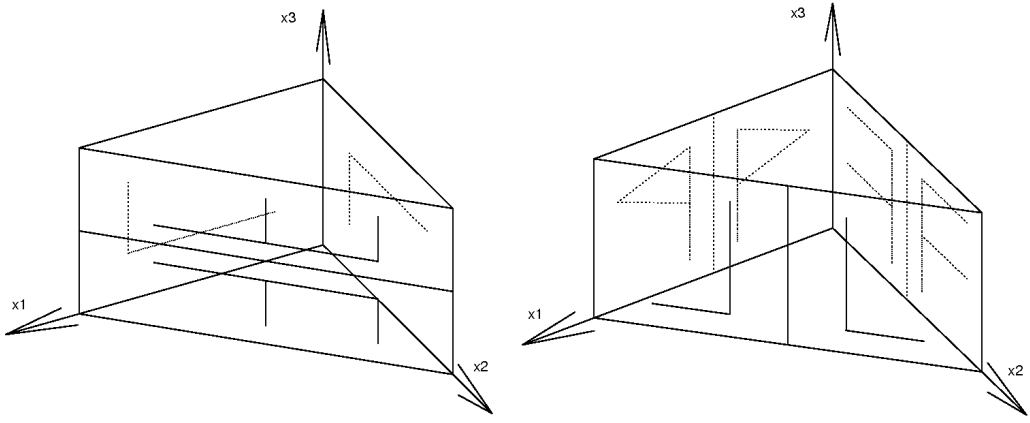}}}
\caption{The underlying spaces of $\O_1$ and $\O_2$ from Example \ref{ex:flat3}}
\label{fig:sunada-indirect}
\end{figure}

No isotropy groups of order greater than two can occur now, since $\tau^2=\rho^2=\Id$,
thus $\Lambda$ is of index two in $\Gamma_1$ and~$\Gamma_2$.
Therefore, all singular points in $\O_1$ and $\O_2$ have isotropy group~$\Z_2$.

The points $x\in\R^3$ which are fixed by an element of the form
$\tau L_\lambda\in\Gamma_1$ must satisfy
$(x_1+x_2,x_1+x_2,2x_3)=-(\lambda_2,\lambda_1,\lambda_3)$.
These are exactly those $x$ with $x_3\in\frac1{2\sqrt2}\Z$ and $x_1+x_2\in\Z$.
Thus (taking identifications into account), the singular set in~$\O_1$ consists
of two copies of $S^1$ of length~$\sqrt2$, corresponding to the horizontal
segments in the face of the fundamental domain denoted by~$\mathsf{F}$ on the left hand side of Figure~\ref{fig:sunada-indirect}.

The orbifold~$\O_2$ is the Riemannian product of a two-dimensional orbifold called
a $4$-pillow or $2222$-orbifold
(here in the form of a square of side length~$\frac12$) and a circle of
length~$\frac1{\sqrt2}$.
Accordingly, its singular set consists of four copies of $S^1$ of length~$\frac1{\sqrt2}$
(corresponding
to the vertical segments over the points $0,\frac12e_1,\frac12e_2,\frac12(e_1+e_2)$; note
that the points $e_1$ and $e_2$ are identified with~$0$).
So, also in this pair of Sunada-isospectral (see below) flat orbifolds,
the singular sets have different topology.

\textit{Sunada isospectrality:}
Define the sublattice $\Lambda':=\spann_\Z\{(1,1,0),(1,-1,0),(0,0,\sqrt2)\}$
of~$\Lambda$, and let $T_{\Lambda'}:=L_{\Lambda'}\backslash\R^3$.
We will see that $\O_1\cong G_1\backslash T_{\Lambda'}$ and
$\O_2\cong G_2\backslash T_{\Lambda'}$
for two eight-element groups $G_1,G_2$ of isometries of~$T_{\Lambda'}$ which are
almost conjugate in a certain finite subgroup of the isometry group of $T_{\Lambda'}$.
Here, we use the symbol~$\cong$ to denote that two orbifolds are isometric.

One easily sees that $\Lambda'$ has index four in~$\Lambda$, and that a 
full set of representatives of $\Lambda/\Lambda'\cong\Z_2\times\Z_2$ is given by
$\{0,e_2,\frac1{\sqrt2}e_3,e_2+\frac1{\sqrt2}e_3\}$.
Since $\Lambda'$ is invariant under $\tau$ and~$\rho$, these isometries of~$\R^3$
descend to isometries $\overline\tau$ and $\overline\rho$ of $T_{\Lambda'}:=L_{\Lambda'}
\backslash\R^3$; trivially, also translations $L_\lambda$ descend to isometries
$\overline{L_\lambda}$ of $T_{\Lambda'}$. Define the groups
$$
G_1:=\{\overline{\tau^jL_\lambda}\mid j\in\{0,1\},\lambda\in\Lambda/\Lambda'\}
\mbox{ and }
G_2:=\{\overline{\rho^jL_\lambda}\mid j\in\{0,1\},\lambda\in\Lambda/\Lambda'\}.
$$
It is not hard to verify that
$$
G_1\backslash T_{\Lambda'}\cong\O_1
\mbox{ and }
G_2\backslash T_{\Lambda'}\cong\O_2.
$$
We are looking for a bijection from $G_1$ to $G_2$ preserving conjugacy classes
in the isometry group of $T_{\Lambda'}$.
Let
$$
A:=\left(\begin{matrix} -1/2&-1/2&1/\sqrt2\\ 1/2&1/2&1/\sqrt2\\ 1/\sqrt2&-1/\sqrt2&0
\end{matrix}\right).
$$
Note that $A\tau A\inv=\rho$ and $A(\Lambda')=\Lambda'$. Let $H\subset O(3)$
be the subgroup generated by $\tau$, $\rho$, and $A$. Note that $H$ is finite
since it preserves the lattice~$\Lambda'$.
Define $\Phi:G_1\to G_2$ by $\Phi(\overline{L_\lambda})
:=\overline{L_\lambda}$ for $\lambda\in\Lambda/\Lambda'$ and
$$
\Phi(\overline\tau):=\overline\rho,\mbox{ }\;
\Phi(\overline{\tau L_{e_2}}):=\overline{\rho L_{\frac1{\sqrt2}e_3}},\mbox{ }\;
\Phi(\overline{\tau L_{\frac1{\sqrt2}e_3}}):=\overline{\rho L_{e_2}},\mbox{ }\;
\Phi(\overline{\tau L_{e_2+\frac1{\sqrt2}e_3}}):=\overline{\rho L_{e_2+\frac1{\sqrt2}e_3}}.
$$
We claim that $\Phi$ preserves conjugacy classes in the finite subgroup
$$
G:=\{\overline{BL_b}\mid B\in H,\; b\in(\Lambda'/4)/\Lambda'\}
$$
of the isometry group of $T_{\Lambda'}$. This follows from the relation
$$
A\tau L_\lambda A\inv=\rho L_{A\lambda}
$$
in connection with the following formulas, where $b:=\frac14e_1-\frac14e_2\in\Lambda^\prime/4$:
\begin{align*}
L_b\inv(\rho L_{Ae_2})L_b&=\rho L_{-\frac1{\sqrt2}e_3}\sim \rho L_{\frac1{\sqrt2}e_3}\\
L_b(\rho L_{A(\frac1{\sqrt2}e_3)})L_b\inv&=\rho L_{e_2}\\
\rho L_{A(e_2+\frac1{\sqrt2}e_3)}&=\rho L_{e_2-\frac1{\sqrt2}e_3}\sim
\rho L_{e_2+\frac1{\sqrt2}e_3}
\end{align*}
Here, the sign $\sim$ between two isometries of~$\R^3$ means that they
differ by a translation in $L_{\Lambda'}$ and thus induce
the same isometry of $T_{\Lambda'}$.
So $\O_1$ and $\O_2$ are indeed Sunada-isospectral;
in particular, they are isospectral on $k$-forms for all~$k$.
\end{example}

\begin{remark}
It is an interesting open question whether there exists a pair of compact
flat orbifolds which are $k$-isospectral for all~$k$ and have different
maximal isotropy orders. Another open question is whether a pair of compact
flat orbifolds which are $k$-isospectral for all~$k$ must necessarily
be Sunada-isospectral.
\end{remark}

\begin{example}
\label{ex:flat4}\
Another pair of Sunada-isospectral orbifolds is given as follows. Let $\Lambda:=2\Z\times2\Z\times2\Z$,
$$
\chi_1:=\left(\begin{matrix} -1&0&0\\0&1&0\\0&0&1\end{matrix}\right),\mbox{ }
\chi_2:=\left(\begin{matrix} -1&0&0\\0&-1&0\\0&0&1\end{matrix}\right),\mbox{ }
\chi_3:=\left(\begin{matrix} 1&0&0\\0&-1&0\\0&0&1\end{matrix}\right),
$$
$$
b_1:=e_1+e_2,\mbox{ }b_2:=0,\mbox{ }b_3:=e_1+e_2,\mbox{ }b^\prime_1:=e_3,\mbox{ }b^\prime_2:=0,\mbox{ }b^\prime_3:=e_3\in\R^3.$$
Set $\rho_j:=\chi_j\circ L_{b_j},\mbox{ }\rho^\prime_j:=\chi_j\circ L_{b^\prime_j}$, $\rho_0=\rho^\prime_0=\Id$ and observe that
$$
\Gamma_1:=\{\rho_j L_\lambda\mid j\in\{0,1,2,3\}, \lambda\in\Lambda\}\mbox{ and }
\Gamma_2:=\{\rho^\prime_j L_\lambda\mid j\in\{0,1,2,3\}, \lambda\in\Lambda\}
$$
are discrete and cocompact subgroups of $O(3)\ltimes\R^3$. Note that the orbifolds
$\O_1:=\Gamma_1\backslash\R^3$ and $\O_2:=\Gamma_2\backslash\R^3$ are not orientable.
For both, a fundamental domain is given by $[0,1]\times[0,1]\times[0,2]$. The
boundary identifications are shown in Figure~\ref{fig:cuboidsJP}, where we omit
the identifications by $L_{2e_3}$ as usual. Note that the underlying topological
space of~$\O_1$
is the product of a projective plane and a circle.

\begin{figure}[!htb]
\centerline{\mbox{\includegraphics*[height=7cm]{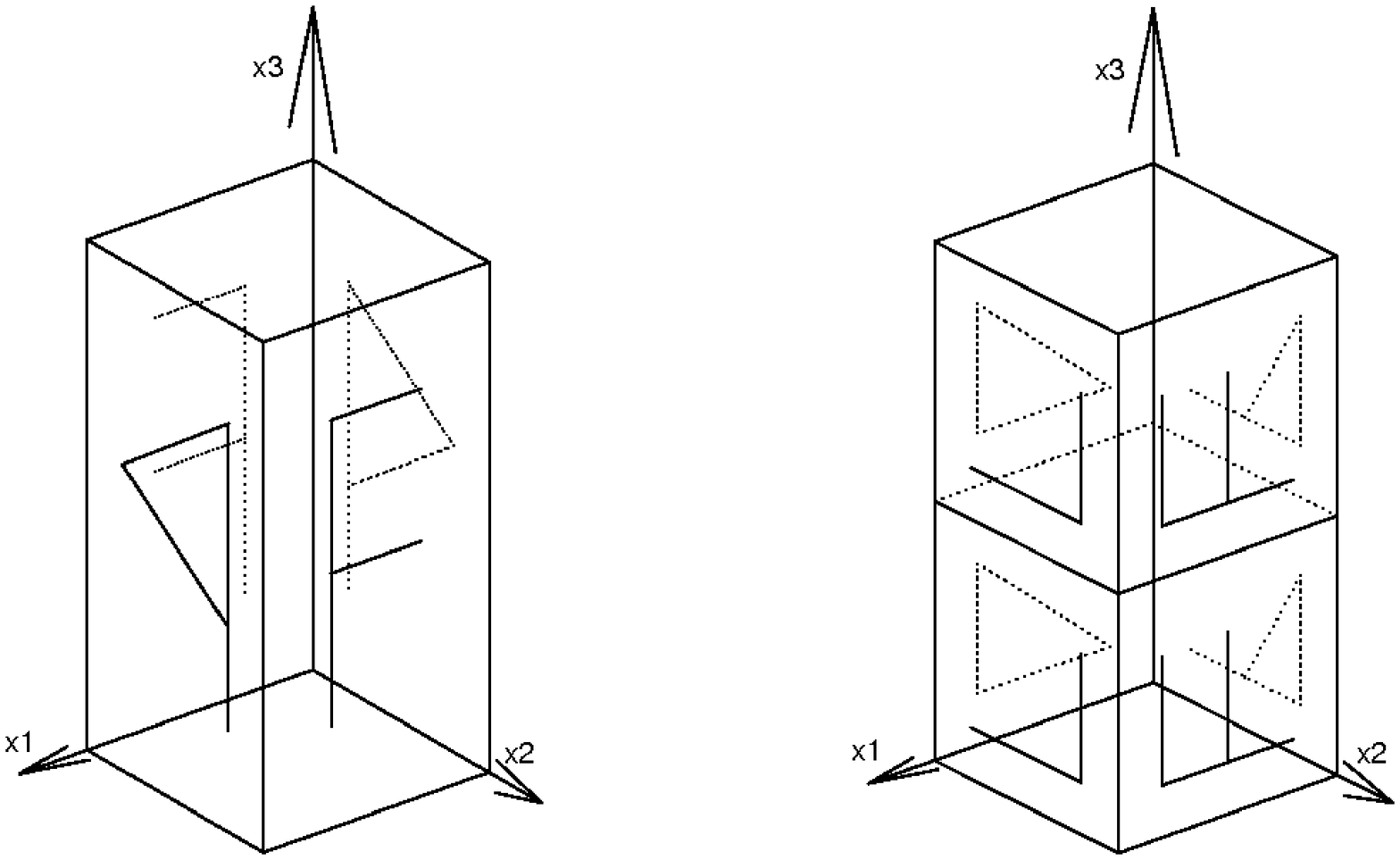}}}
\caption{The underlying spaces of $\O_1$ and $\O_2$ from Example \ref{ex:flat4}}
\label{fig:cuboidsJP}
\end{figure}

Using the notation introduced at the beginning of this section, we note that $\O_i\cong \bar{\Gamma}_i\backslash T_\Lambda$ for $i=1,2$, where $\bar\Gamma_1$, $\bar\Gamma_2$ are
the following subgroups of the isometry group of $T_\Lambda$:
$$\bar{\Gamma}_1=\{\overline{\rho_j}\mid j\in\{0,1,2,3\}\},\mbox{ }\bar{\Gamma}_2=\{\overline{\rho^\prime_j}\mid j\in\{0,1,2,3\}\}.$$
It is not difficult to see that the groups $\bar{\Gamma}_1$ and $\bar{\Gamma}_2$ are almost
conjugate in the finite group~$G$ generated by $\bar{\Gamma}_1$,
$\bar{\Gamma}_2$ and $\{\overline{B L_b}\mid B\in P(3),b\in (\Lambda/4)/\Lambda\}$,
where $P(3)\subset O(3)$ denotes the group of permutation matrices. Hence, $\O_1$
and $\O_2$ are Sunada-isospectral. Alternatively, one can apply the methods developed in \cite{MR:2002}, Section 3, to verify that the two orbifolds are Sunada-isospectral.

However, $\O_1$ and $\O_2$ are not isometric; in fact,
their respective singular sets have different numbers of components.
For each $i=1,2$ the points in $\R^3$ which are fixed by nontrivial elements
of $\Gamma_i$ are given by the set $\Z\times\Z\times\R$. Each of these points
is fixed by exactly one nontrivial group element and thus has isotropy $\Z_2$.
Taking identifications into account (recall Figure~\ref{fig:cuboidsJP}), we
observe that in $\O_1$ the singular set consists of two copies of $S^1$ of length two, whereas in
$\O_2$ it consists of four copies of $S^1$ of length one.
\end{example}

Finally, we present another pair of (non-Sunada) isospectral orbifolds with properties
similar to the pair from Example~\ref{ex:flat2}, this time with nonisomorphic maximal
isotropy groups of order six.

\begin{example}
\label{ex:flat5}\
Let $\Lambda:=\spann_\Z\{(2,0,0),(1,\sqrt{3},0),(0,0,1)\}$ and
$$
H:=\left(\begin{matrix} 1/2&-\sqrt{3}/2&0\\\sqrt{3}/2&1/2&0\\0&0&1\end{matrix}\right),\mbox{ }
R:=\left(\begin{matrix} 1&0&0\\0&-1&0\\0&0&-1\end{matrix}\right).
$$
Note that $H$ is just the rotation by $\pi/3$ around the $x_3$-axis. Now
$$
\Gamma_1:=\{H^j L_\lambda\mid j\in\{0,\ldots,5\}, \lambda\in\Lambda\},\mbox{ }
\Gamma_2:=\{H^{2j} R^k L_\lambda\mid j\in\{0,1,2\}, k\in\{0,1\}, \lambda\in\Lambda\}
$$
are crystallographic groups acting on $\R^3$. For both $i=1,2$, a fundamental domain
of the action of~$\Gamma_i$ on $\R^3$ is given by the prism of height one over the triangle with
vertices $(0,0,0)$, $(2,0,0)$, $(1,1/\sqrt{3},0)$ (compare Figure~\ref{fig:hexorbis} where we
again omit the identifications by $L_{e_3}$).

\begin{figure}[!htb]
\centerline{\mbox{\includegraphics*[height=7cm]{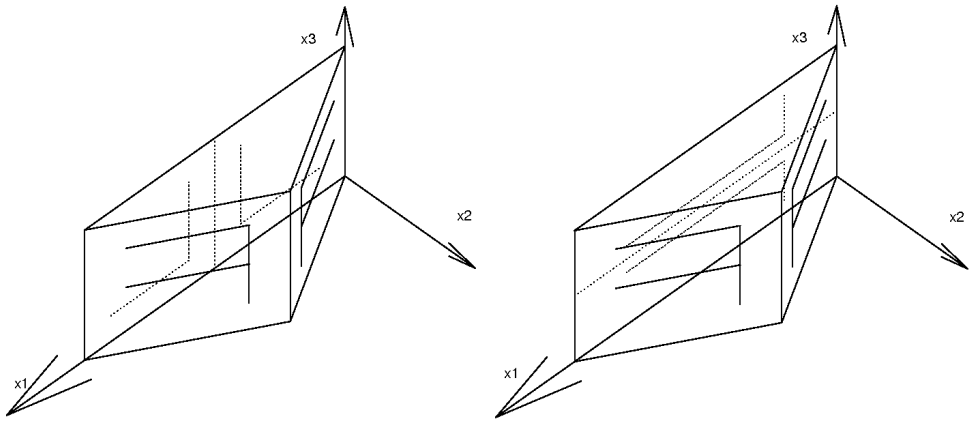}}}
\caption{The underlying spaces of $\O_1$ and $\O_2$ from Example \ref{ex:flat5}}
\label{fig:hexorbis}
\end{figure}

Using Theorem~\ref{thm:dimformula} one
shows that the two orbifolds $\O_1:=\Gamma_1\backslash\R^3$
and $\O_2:=\Gamma_2\backslash\R^3$ are isospectral on functions but not on $1$-forms.
It is not hard to verify that the maximal isotropy group is $\Z_6$ in the case of $\O_1$
and $D_6$ (the dihedral group with six elements) in the case of $\O_2$. Just as in
Example~\ref{ex:flat2}, the sets of points with maximal isotropy have different dimensions:
In~$\O_1$, it is a circle of length one (the image of the $x_3$-axis), while in~$\O_2$
it consists of only two points (the images of $(0,0,0)$ and $(0,0,1/2)$).
Note that $\O_1$ is the product of a $236$-orbifold with a circle of length one.
So its other nontrivial isotropy groups are $\Z_2$ and $\Z_3$, and the corresponding
singular points each time form another circle of length one. In $\O_2$ there
are two open segments of length two consisting of points with isotropy group~$\Z_2$ (corresponding
to the horizontal segments in Figure~\ref{fig:hexorbis}). The set of points with
isotropy $\Z_3$ consists of the open segment of length $1/2$
which joins the two points with maximal isotropy and of the circle of length one
corresponding to the vertical edge through the point $(1,1/\sqrt3,0)$.

\end{example}

\end{document}